# Graceful Labeling For bipartite graceful Graphs and related Graphs☆

Chunfeng Liu[a], XiuyingWang[b] ,Guozhu Liu[b] *

([a]College of Science, Liaoning University of Technology, Liaoning,Jinzhou，121013)

( [b]College of Information Science and Technology，Qingdao University of Science and Technology，Shandong, Qingdao 266000)

**Abstract:** The concept of graceful labels was proposed by Rosa, scholars began to study graceful labels of various graphs and obtained relevant results.In this paper，let $G$ is a bipartite graceful graph, we proved that $S(t \bullet G)$, $P_n^t(tn \bullet G)$, $C(t \bullet G)$ and $G^\star$ are graceful labelings.

## 1. INTRODUCTION AND LEMMEL

The graceful labeling was introduced by A Rosa [1] in 1967. Golomb [2] proved that the complete bipartite graph is graceful. Barrientos [3] proved that union of complete bipartite graphs is also graceful. Vaidya [4] introduced a star of cycle. Kaneria and Makadia [5] proved that star of a cycle is graceful. Kaneria [6] proved that join sum of path union and star are graceful graphs. The detail survey of graph labeling refer Gallian [7].

Let $G$ is a bipartite graceful graph, In this paper, we proved that $S(t \bullet G)$, $P_n^t(tn \bullet G)$, $C(t \bullet G)$ and $G^\star$ are graceful graphs. We assume $G$ is a simple undirected, finite graph, with vertices and edges. For all terminology and notations we follow Harary [8]. We shall give brief summary of definitions which are useful in this paper.

**Definition1.1** A function $f$ is called graceful labeling of a $(p, q)$-graph $G$ if $f$: $V(G) \to \{0,1,2,\cdots, q\}$ is injective and induced function $f^*$: $E(G) \to \{0,1,2,\cdots, q\}$ defined as $f(e)=|f(u)-f(v)|$ is bijective for every edge $e=uv$. A graph $G$ is called graceful graph if it admits a graceful labeling.

**Definition1.2** Let $(p,q)$-graph $G$ is a bipartite graph with the bipartition $(X,Y)$, $f$ is a graceful labeling of $G$ , if $max\{f(x)|x \in X\}<min\{f(y)|y \in Y\}$,then $f$ called is a bipartite graceful labeling of $G$. if $G$ admits a bipartite graceful labeling, then call $G$ is a bipartite graceful graph.

**Definition1.3** Let $G_1, G_2,\cdots, G_n$ ($n\geq 2$) be $n$ graph ,Then the graph obtained by adding an edge from $G_i$ to $G_{i+1}$($1\leq i\leq n-1$) is called path union of graph $G_1, G_2,\cdots, G_n$. We shall denote such graph by $P(G_1, G_2,\cdots, G_n)$.

If $G_1=G_2=\cdots=G_n=G$, we shall denote by $P(nt \bullet G)$.

**Definition1.4**[10] A graph obtained by replacing each vertex of $K_{1,t}$ except the apex vertex by the graphs $G_1, G_2,\cdots, G_t$ is known as open star of graphs. We shall denote such graph by $S(G_1, G_2,\cdots, G_t)$.

If we replace each vertices of $K_{1,n}$ except the apex vertex by a graph $G$. i.e. $G_1=G_2=\cdots=G_t=G$, such open star of a graph, we shall denote by $S(t \bullet G)$.

**Definition1.5**[10] A graph is obtained by replacing each edge of $K_{1,t}$ by a path of $P_n$ length $n$ on $n+1$ verticesis called one point union for $t$ copies of path $P_n$. We shall denote such graph $G$ by $P_n^t$.

**Definition 1.6**[10] A graph obtained by replacing each vertices of except the central vertex by the graphs $G_1, G_2,\cdots, G_{tn}$ is known as one point union for path of graphs. We shall denote such graph $G$ by $P_n^t(G_1, G_2,\cdots, G_{tn})$,where $P_n^t$ is the one point union of $t$ copies of path $P_n$. If we replace each vertices of $P_n^t$ except the central vertex by a graph$G$, i.e. $G=G_1=G_2=\cdots=G_{tn}$ such one point union of path graph, we shall denote it by $P_n^t(tn,G)$.

**Definition1.7**[11] For a cycle $C_n$, each vertices of $C_n$ is replace by connected graphs $G_1, G_2,\cdots, G_n$ is known as *cycle of graphs* and we shall denote it by $C(G_1, G_2,\cdots, G_n)$. If we replace each vertices by graph $G$,i.e. $G=G_1=G_2=\cdots=G_n$, such cycle of a graph $G$, we shall denote it by $C(n \bullet G)$.

☆This work was Supported by the National Natural Science Foundation of China (61773107,61603168)

*Corresponding Author at College of Information Science and Technology，Qingdao University of Science and Technology，Shandong, Qingdao 266000, E-mail：lgz_0228@163.com

**Definition1.8**[6] A graph obtained by replacing each vertex of star $K_{1,n}$ by a graph $G$ of $n$ vertices is called star of $G$ and it is Denoted by $G^{\star}$. The graph $G$ which replaced at the center of $K_{1,n}$ we call the central copy of $G^{\star}$.

Let $G$ is a bipartite graceful graph, In this paper, We also proved that $P(n\bullet G)$,$S(t\bullet G)$ $P_n^{t}(tn\bullet G)$, $C(t\bullet G)$ and $G^{\star}$ are graceful graphs.

Obviously, there are following conclusions:

**Lemma 1.1** $P_n$ is bipartite graceful graph.

**Lemma 1.2** $C_n$ ( $n\equiv 0\ (mod\ 4)$) is bipartite graceful graph.

**Lemma 1.3** $K_{m,n}$ ( $m,n\in N$) is bipartite graceful graph.

**Lemma 1.4** $P_m\times P_n$ ( $m,n\in N$) is bipartite graceful graph.

## 2 MAIN RESULTS

**Theorem 2.1** Let $G$ is bipartite graceful graph and $G_i(i=1,2,\cdots,t)$ is be $t$ copies of graph $G$, Then $P(n\bullet G)$ is graceful.

**Proof:** Let $(p_0,q_0)$-graph $G$ is a bipartite graph with the bipartition $(V_1,V_2)$ ,and $V_1=\{u_i|i\in[1,m]\}$, $V_2=\{v_i|j\in[1,r]\}$, $f_0$ is a bipartite graceful labeling of $G$ , $f_0(u_i)<f_0(u_{i+1})$, $i\in[1,m-1]$, $f_0(v_j)<f_0(v_{j+1})$, $j\in[1,r-1]$, and $f_0(u_m)<f_0(v_1)$ . Let $u_{l,i}(i\in[1,m])$, $v_{l,j}(j\in[1,r])$ be the vertices of $G$ which is $l^{th}$ copy of the path union of $n$ copies of $G$. We shall join $v_{l,1}$ with $u_{l,1}$.

We shall define labeling function of $P(n\bullet G)$ $f$ as follows:

$f$: $P(n\bullet G)\rightarrow[0,q]$ ,Where $q=n(q_0+1)-1$ defined by,

$f(u_{1,i})=f_0(u_i)$, $\quad i\in[1,m]$;

$f(v_{1,j})=q-q_0+f_0(v_j)$, $\quad j\in[1,r]$;

$f(u_{l,i})=f(u_{l-1,m})+f_0(u_i)+1$, $\quad i\in[1,m]$, $l\in[2,n]$;

$f(v_{l,j})=f(v_{l-1,1})+f_0(v_j)-q_0-1$, $\quad j\in[1,r]$, $l\in[2,n]$.

Above labeling pattern give rise graceful labeling to the path union graph $P(n\bullet G)$ of $n$ copies of $G$.

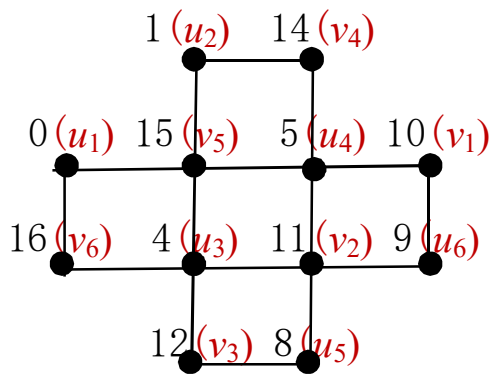

Figure 1.Graph $P_{l4}$ and its bipartite graceful labeling

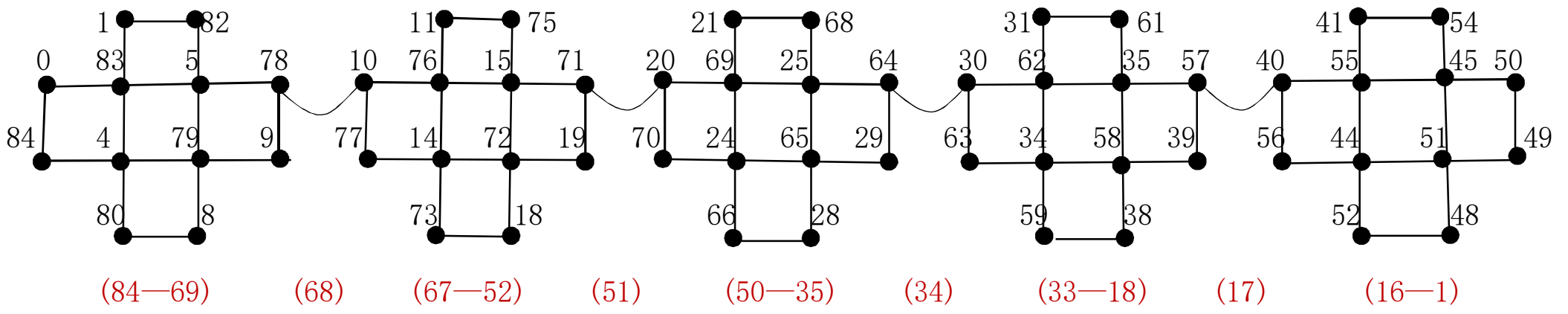

Figure 2.Graph $P(5\bullet P_{l4})$ and its bipartite graceful labeling

By theorem2.1 and lemma 1.1−1.4, respectively,have following consequences.

**Corollary 2.1.1** $P(t\bullet P_n)$ is graceful.

**Corollary 2.1.2**[13] $P(t\bullet C_n)$ is graceful graph, where $n\equiv 0\ (mod\ 4)$.

**Corollary 2.1.3** $P(t\bullet K_{m,n})$ is graceful.

**Corollary 2.1.4**[12] $P(t\bullet P_m\times P_{,n})$ is graceful.

**Theorem 2.2** Let $G$ is bipartite graceful graph and $G_i(i=1,2,\cdots,t)$ is be $t$ copies of graph $G$, Then $S(t\bullet G)$ is graceful.

**Proof:** Let $(p_0,q_0)$-graph $G$ is a bipartite graph with the bipartition $(V_1,V_2)$ ,and $V_1=\{u_i\,|\,i\in[1,m]\}$, $V_2=\{v_i\,|\,j\in[1,n]\}$, $f_0$ is a bipartite graceful labeling of $G$ , $f_0(u_i)<f_0(u_{i+1})$, $i\in[1,m-1]$, $f_0(v_j)<f_0(v_{j+1})$, $j\in[1,n-1]$, and $f_0(u_m)<f_0(v_1)$. Let $S(t\bullet G)$ be a graph obtained by replacing each vertices of $K_{1,t}$ except the apex vertex by the graph $G$.Let $v_0$ is the apex vertex of $K_{1,t}$ . i.e. it is central vertex of the graph $S(t\bullet G)$. Let $u_{li}(1\le i\le m)$, $v_{lj}(1\le j\le n)$ be the vertices of $l^{th}$ copy of $G^{(l)}$ of $G$ in $S(t\bullet G),\forall l\in[1,t]$.

We shall join $v_{lm}$ with the vertex $v_0$ by an edge to form the open star of graphs $S(t\bullet G)$, $\forall l\in[1,t]$. where $q=t(q_0+1)$ . We define labeling function $f$: $V(G)\rightarrow[1,\ q]$ as follows:

$f(v_0)=0$;

$f(u_{1,i})=f_0(u_i)+1,\qquad i\in[1,m]$;

$f(v_{1,j})=q-q_0+f_0(v_j),\qquad j\in[1,n]$;

$f(u_{2,i})=f(u_{1,i})+q-q_0-1,\qquad i\in[1,m]$;

$f(v_{2,j})=f(v_{1,j})-q+q_0+1,\qquad j\in[1,n]$;

$f(u_{i,j})=f(u_{l-2,j})-(-1)^l(q_0+1),\ i\in[1,m],\ l\in[3,t]$;

$f(v_{i,j})=f(v_{l-2,j})+(-1)^l(q_0+1),\ j\in[1,n],\ l\in[3,t]$

Above labeling pattern give rise graceful labeling to the graph and so it is a graceful graph.

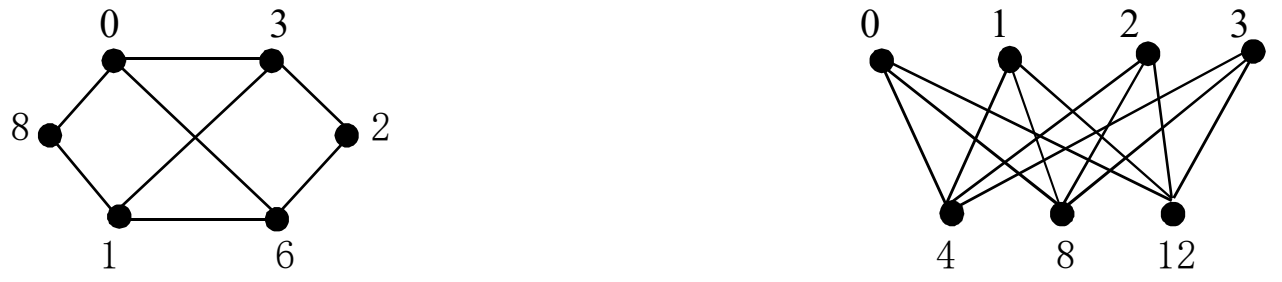

(a) Graph Q and its bipartite graceful labeling (b) Graph $K_{4,3}$ and its bipartite graceful labeling

Figure 3

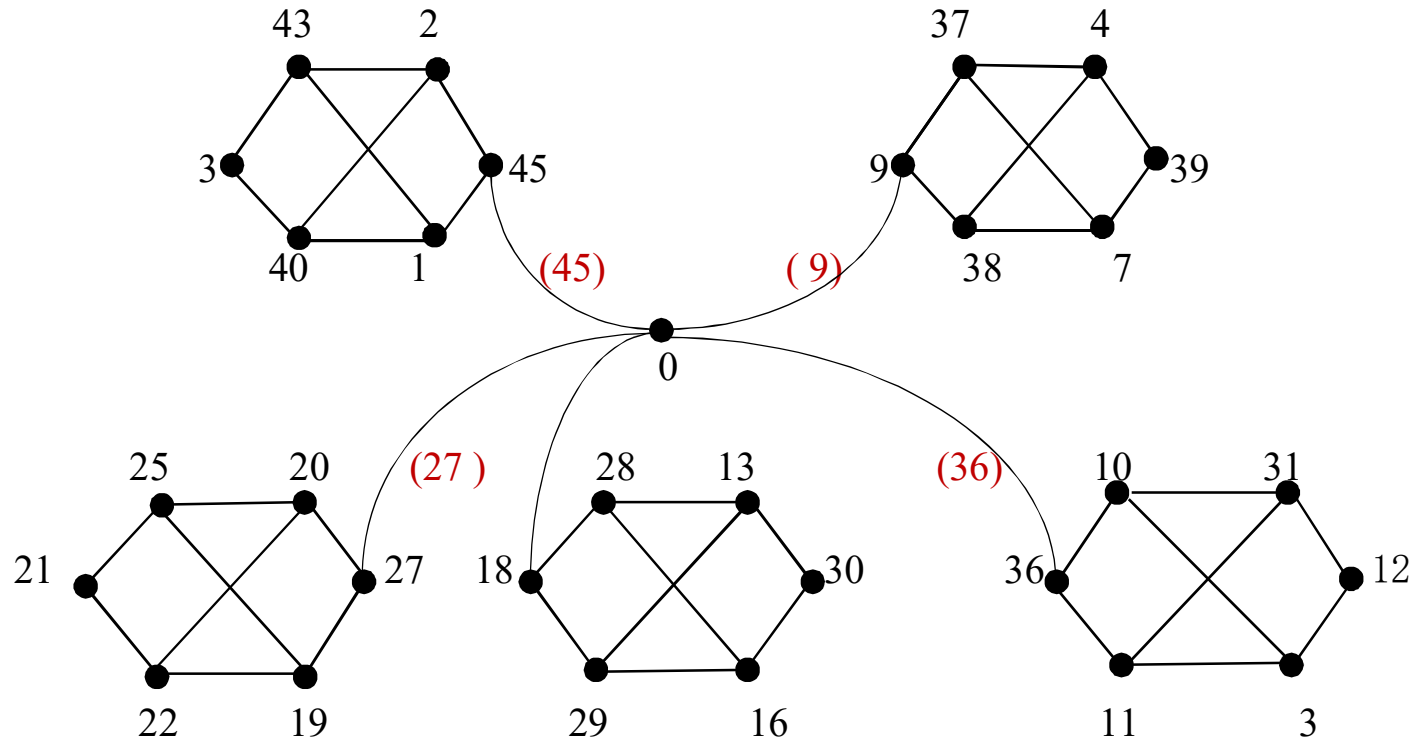

Figure 4. A graph obtained by open star of $Q$ and its graceful labeling .

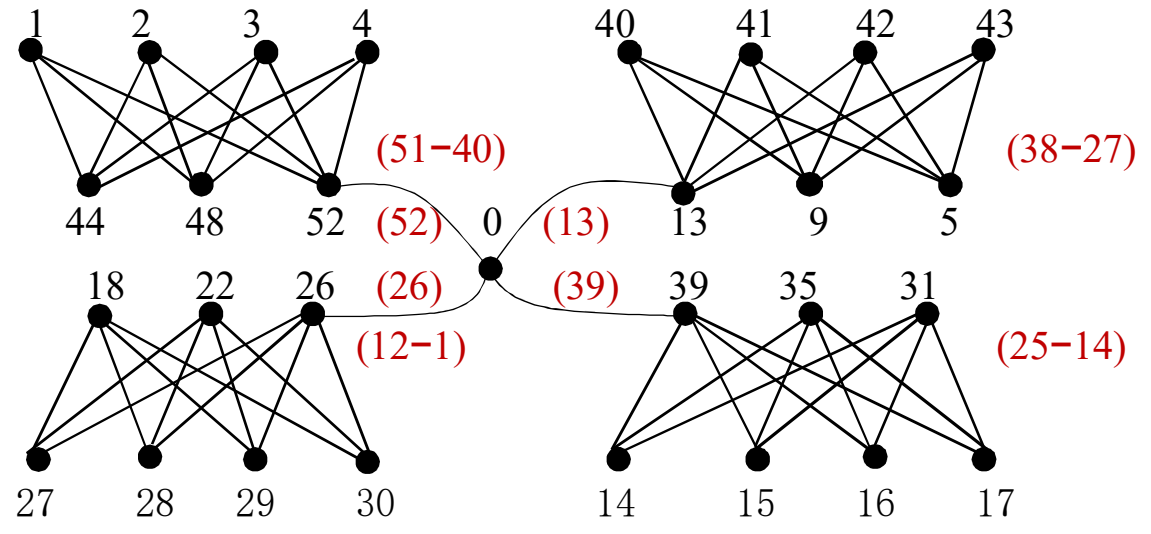

Figure 5. A graph obtained by open star of $K_{4,3}$ and its graceful labeling .

By theorem2.2 and lemma 1.1−1.4, there are following consequences respectively.

**Corollary 2.2.1** $S(t\cdot P_n)$ is graceful.

**Corollary 2.2.2**[10] $S(t\cdot C_n)$ is graceful graph, where $n\equiv 0\ (mod\ 4)$.

**Corollary 2.2.3**[9] $S(t\cdot K_{m,n})$ is graceful.

**Corollary 2.2.4**[9] $S(t\cdot P_m\times P_n)$ is graceful.

**Theorem2.3** Let $G$ is bipartite graceful graph and $G_i$ ($i$= 1,2,⋯, $t$) is $ith$ copies of graph $G$, Then $P_n^t(tn\cdot G)$ is graceful.

**Proof** Let $(p_0,q_0)$-graph $G$ is a bipartite graph with the bipartition $(V_1,V_2)$ ,and $V_1=\{u_i|i\in[1,m]\}$, $V_2=\{v_i|j\in[1,r]\}$, $f_0$ is a bipartite graceful labeling of $G$ , $f_0(u_i)<f_0(u_{i+1})$, $i\in[1,m-1]$, $f_0(v_j)<f_0(v_{j+1})$, $j\in[1,r-1]$, and $f_0(u_m)<f_0(v_1)$ .Let be a graph obtained by replacing each vertices of $P_n^t$ except the central vertex by the graph $G$ .i.e. $P_n^t(tn\cdot G)$ is the graph obtained by replacing each vertices of $K_{1,t}$ except the apex vertex by the path union of $n$ copies of the graph $G$. Let $u_0$ be the central vertex for the graph $P_n^t(tn\cdot G)$. Let $u_{s,l,i}(i\in[1,m])$, $v_{s,l,j}(j\in[1,r])$ be the vertices of $G$ which is $l^{th}$ copy of the path union of $n$ copies of $G$ lies in $s^{th}$ branch of the graph $P_n^t(tn\cdot G)$,$l\in[1,n]$ and $s\in[1,t]$.

First we shall join $v_{s,l,m}$ with the vertex $u_{s,l+1,1}$ by an edge to form the path union of $n$ copies of $G$ for $s^{th}$ branch of $P_n^t(tn\cdot G)$, $l\in[1,n-1]$and $s\in[1,t]$.Now we shall join $u_{s,1,1}$ with the vertex $u_0$ by an edge to form the one point union for path of graphs $P_n^t(tn\cdot G)$ , $s\in[1,t]$.We shall define labeling function $f$ for the first copy(branch) of the path union of $n$ copies of $G$ as follows:

$f$: $P(n\cdot G)\rightarrow[0,q]$ ,Where $q=n(q_0+1)-1$ defined by,

$f(u_{1,1,i})=f_0(u_i)$, $\quad i\in[1,m]$;

$f(v_{1,1,j})=q-q_0+f_0(v_j)$, $\quad j\in[1,r]$;

$f(u_{1,l,i})=f(u_{1,l\text{-}1,m})+f_0(u_i)+1$, $\quad i\in[1,m]$, $l\in[2,n]$;

$f(v_{1,l,j})=f(v_{1,l\text{-}1,1})+f_0(v_j)-q_0-1$, $\quad j\in[1,r]$, $l\in[2,n]$.

Above labeling pattern give rise graceful labeling to the path union of $n$ copies of $G$, which lies in first branch of $P_n^t(tn\cdot G)$. Now we shall define labeling $g$:$V(P_n^t(tn\cdot G))\rightarrow[0,Q]$, where $Q=|E(P_n^t(tn\cdot G)|=tn(q_0+1)$ as follows:

$g(u_0)=0$;

$g(u_{1,l,i})=f(u_{1,l,i})+1$, $\quad i\in[1,m]$, $l\in[1,n]$;

$g(v_{1,l,j})=f(v_{1,l,j})+Q-q$, $\quad j\in[1,r]$, $l\in[1,n]$;

$g(u_{2,l,i})=g(u_{1,l,i})+Q-q-1$, $\quad i\in[1,m]$, $l\in[1,n]$;

$g(v_{2,l,j})=g(v_{1,l,j})-Q+q+1$, $\quad j\in[1,r]$, $l\in[1,n]$;

$g(u_{s,l,i})=g(u_{s\text{-}2,l,i})-(-1)^s(q+1)$, $i\in[1,m]$, $l\in[1,n]$, $s\in[3,t]$;

$g(v_{s,l,j})=g(v_{s\text{-}2,l,j})+(-1)^s(q+1)$, $j\in[1,r]$, $l\in[1,n]$, $s\in[3,t]$.

Above labeling pattern give rise graceful labeling to the graph $P_n^t(tn\cdot G)$ and so it is a graceful graph.

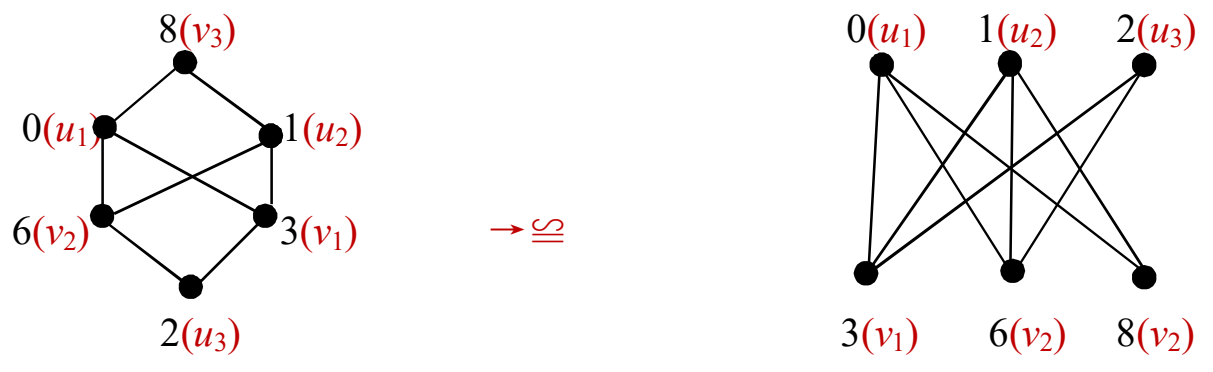

Figure 6 . A graph $Q$ and its bipartite graceful labeling

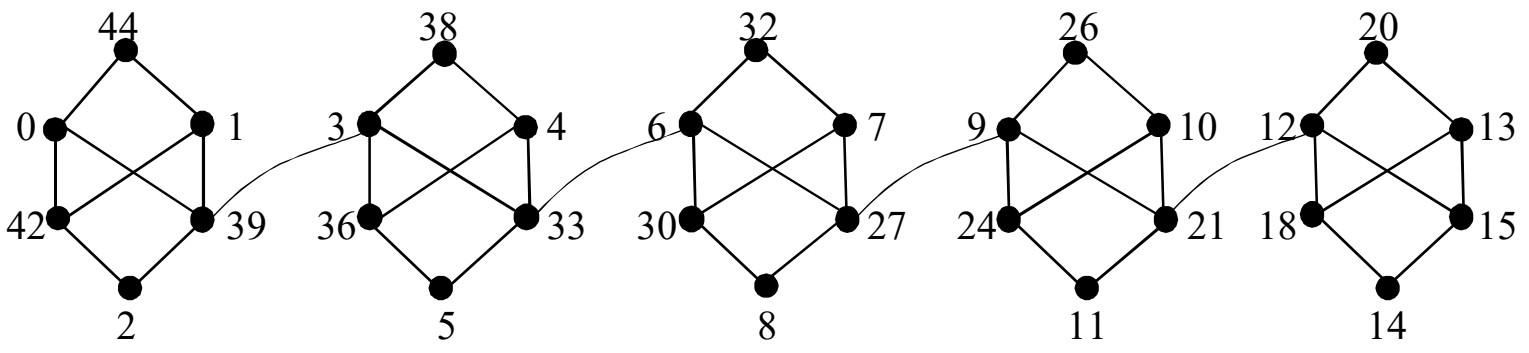

Figure 7 . The path union 5 copies of $Q$ and its graceful labeling

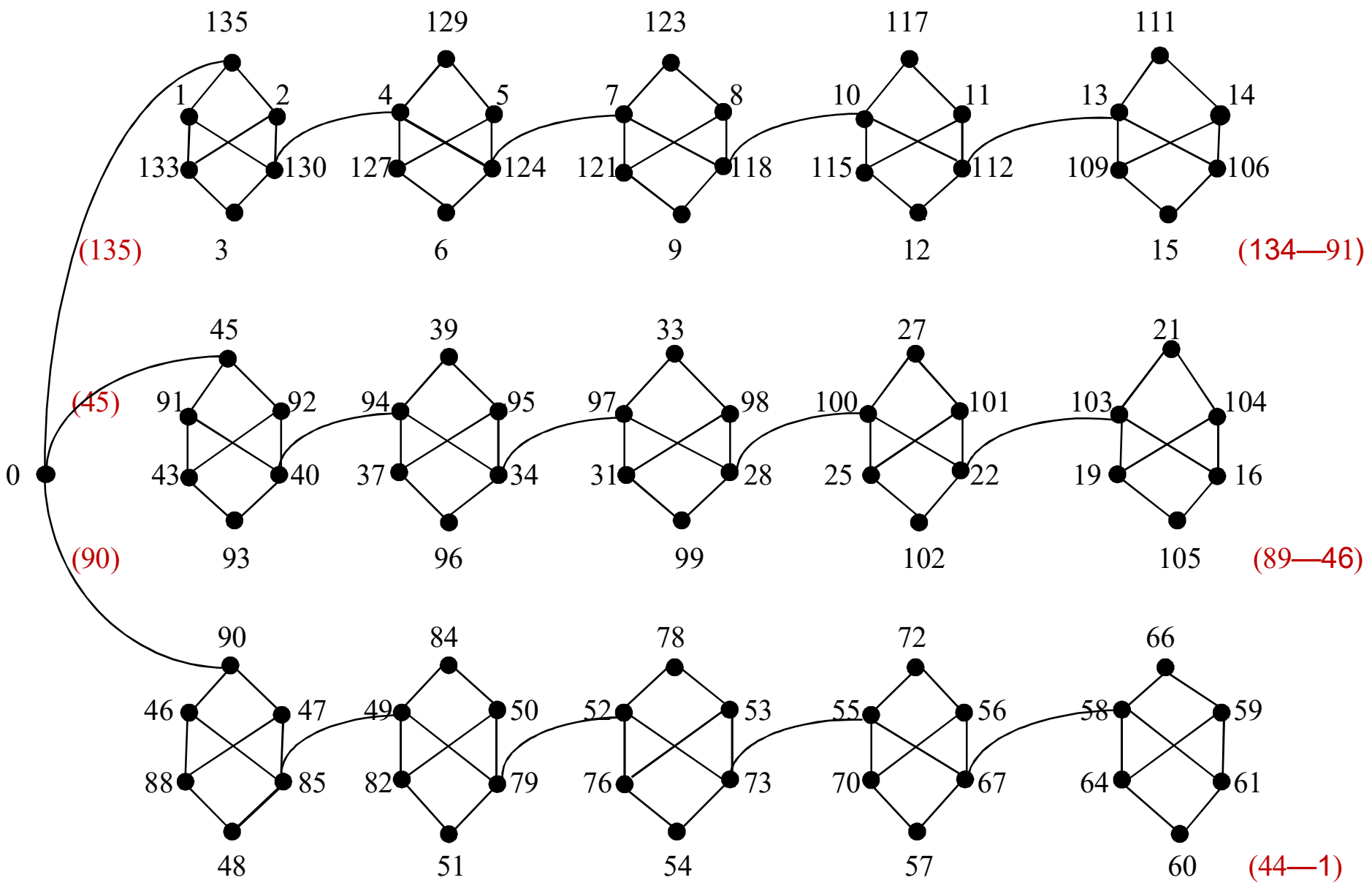

Figure 8 . $P_5^3(3\times 5 \cdot Q)$ and its graceful labeling

By theorem2.3 and lemma 1.1−1.4,respectively,have following consequences.

**Corollary 2.3.1** $P_n^t(tn \cdot P_n)$ is graceful.

**Corollary 2.3.2**[10] $P_n^t(tn \cdot K_{1,m})$ is graceful.

**Corollary 2.3.3**[10] $P_n^t\ (tn \cdot C_m)$ is graceful, where $m\equiv 0(mod4)$.

**Corollary 2.3.4**[9] $P_n^t(tn \cdot K_{m,r})$ is graceful.

**Corollary 2.3.5**[10] $P_n^t(tn \cdot P_r\times P_s)$ is graceful.

**Theorem 2.4** Cycle $C(t \cdot G)$ ($t\equiv 0(mod2)$)of bipartite graceful graph $G$ is graceful graphs.

**Proof** Let Let $(p_0,q_0)$-graph $G$ is a bipartite graph with the bipartition $(V_1,V_2)$ ,and $V_1=\{u_i|i\in[1,m]\}$, $V_2=\{v_i|j\in[1,n]\}$, $f_0$ is a bipartite graceful labeling of $G$ , $f_0(u_i)<f_0(u_{i+1})$, $i\in[1,m-1]$, $f_0(v_j)<f_0(v_{j+1})$, $j\in[1,n-1]$, and $f_0(u_m)<f_0(v_1)$ .$C(t \cdot G)$ be a graph which contains $t$ copies of the bipartite graceful graph $G$,where $t\equiv 0(mod2)$.Let $u_{i,j}(i\in[1,m])$ and $v_{i,j}(j\in[1,n])$ be vertices of $i^{th}$ copy of $G$, $i\in[1,t]$. Now join $v_{i,1}$ with $u_{i+1,1}(i\in[1,(t/2)-1])$, join $v_{t/2,1}$ with $v_{(t/2)+1,1}$, join $u_{i,1}$ with $v_{i+1,1}(i\in[(t/2)+1,t-1])$ and join $u_{t,1}$ with $u_{1,1}$ by an edge to form cycle of graphs $C(t \cdot G)$.

We define the labeling function $f:V(C(tG))\rightarrow[0,q]$, where $q=t(q_0+1)$ as follows:

$f(u_{1,j})=f_0(u_j)$, $j\in[1,m]$;

$f(v_{1,j})=q-f_0(v_n)+f_0(v_j)$, $j\in[1,n]$;

$f(u_{l,i})=f(u_{l-1,m})+1+f_0(u_j)$, $j\in[1,m]$, $l\in[2, t/2]$;

$f(v_{l,j})=f(v_{l-1,1})-1-f_0(v_n)+f_0(v_j)$, $j\in[1,n]$, $l\in[2, t/2]$;

$f(u_{(t/2)+1,i})=f(v_{t/2,1})-2-f_0(u_m)+f_0(u_j)$, $j\in[1,m]$;

$f(v_{(t/2)+1,j})=f(u_{t/2,m})+1-f_0(v_1)+f_0(v_j)$,   $j\in[1,n]$;

$f(u_{l,i})=f(u_{l-1,1})-1-f_0(u_m)+f_0(u_j)$,   $j\in[1,m]$, $l\in[(t/2)+2, t]$;

$f(v_{l,j})=f(v_{l-1,n})+1-f_0(v_1)+f_0(v_j)$,   $j\in[1,n]$, $l\in[(t/2)+2, t]$.

Above labeling pattern give rise graceful labeling to the graph $C(tG)$ and so it is a graceful graph.

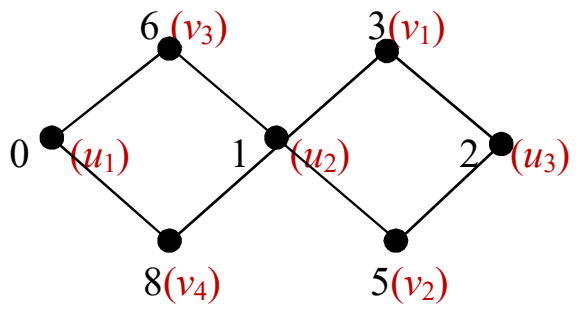

Figure 9 . $G_0$ and its bipartite graceful labeling

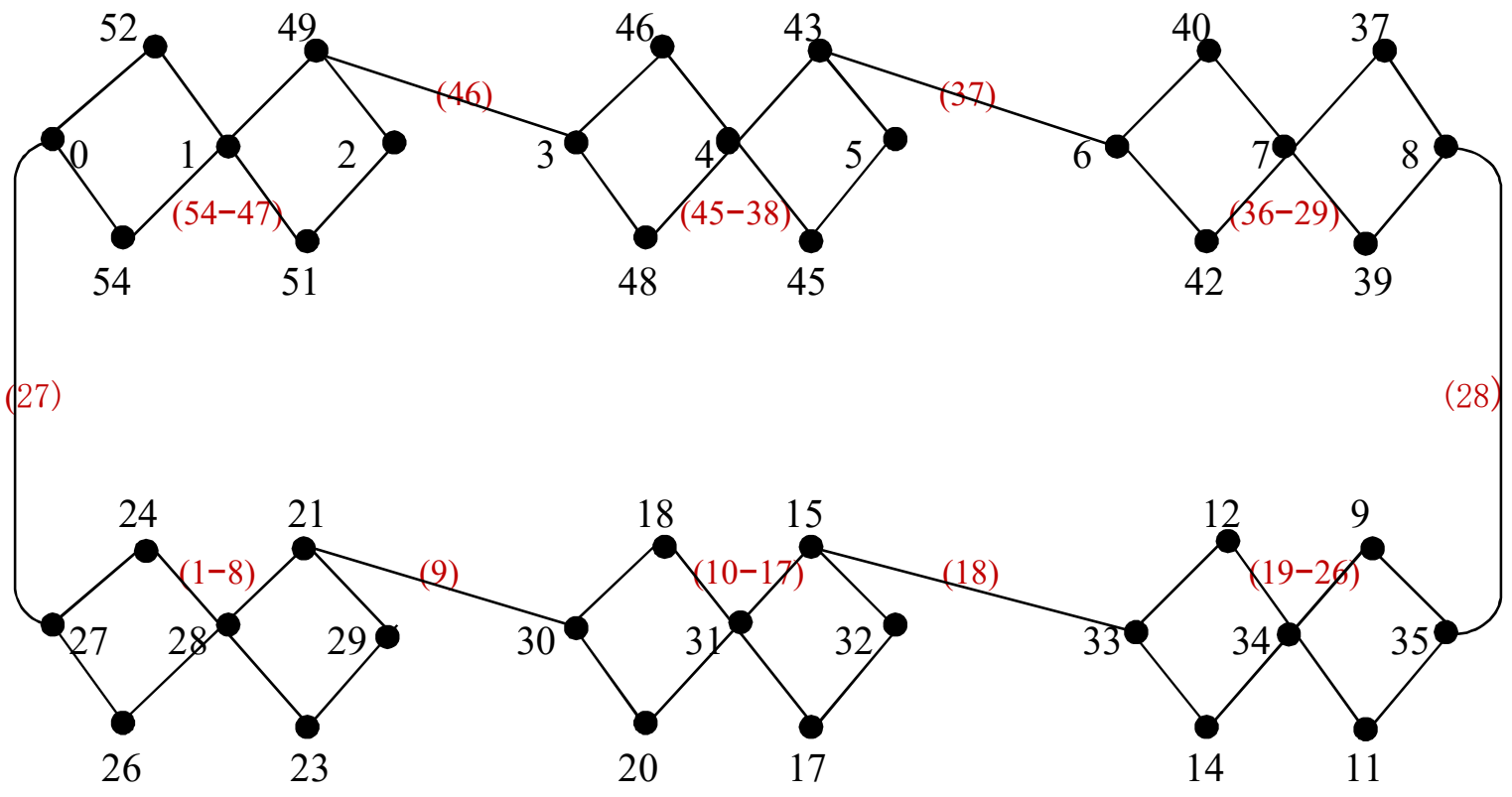

Figure 10. Cycle $C(6 \cdot G_0)$ and its graceful labeling

By theorem2.4 and lemma 1.1−1.4,respectively,have following consequences.

**Corollary2.4.1** Cycle of path graph $C(t \cdot P_n)$ is graceful.

**Corollary2.4.2**[14] Cycle of $C_{4m}$ graphs $C(t \cdot C_{4m})$ , $t\equiv 0(mod2)$ is graceful graphs.

**Corollary2.4.3**[11] Cycle of complete bipartite graphs $C(t \cdot K_{m,n})$, $t\equiv 0(mod2), m,n\in N$ is graceful graphs.

**Corollary2.4.5**[12] Cycle of grid graph $C(t \cdot P_r\times P_s)$ is graceful.

**Theorem2.5**. $G^{\star}$ the star of a bipartite graceful graph $G$, is graceful graphs.

**Proof** Let Let $(p_0,q_0)$-graph $G$ is a bipartite graph with the bipartition $(V_1,V_2)$ ,and $V_1=\{u_i|i\in[1,m]\}$, $V_2=\{v_i|j\in[1,n]\}$, $f_0$ is a bipartite graceful labeling of $G$ , $f_0(u_i)<f_0(u_{i+1})$, $i\in[1,m-1]$, $f_0(v_j)<f_0(v_{j+1})$, $j\in[1,n-1]$, and $f_0(u_m)<f_0(v_1)$ . Let $u_{0,i}$ $(1\leqslant i\leqslant m)$, $v_{0,j}$ $(1\leqslant j\leqslant n)$ be vertices of central copy of $G$ and $u_{l,i}$ $(1\leqslant i\leqslant m)$, $v_{l,j}$ $(1\leqslant j\leqslant n)$ be vertices of other copies of $G$,$\forall$ 1 = 1,2, ⋯,$m+n$.We define labeling function $f: V\rightarrow\{0,1,\cdots,q\}$, where $q=(m+n+1)q_0+m+n$ as follows:

$f(u_{0,i})=f_0(u_i)$,   $i\in[1,m]$;

$f(v_{0,j})=q-f_0(v_n)+f_0(v_j)$   $j\in[1,n]$;

$f(u_{1,i})=f(v_{0,1})-1-f_0(u_m)+f_0(u_i)$,   $i\in[1,m]$;

$f(v_{1,j})=f_0(v_j)$,   $j\in[1,n]$;

$f(u_{l,i})=f(u_{l-2,i})+q_0+1$,   $i\in[1,m]$, $l\in[1,m+n]$ and $l$ *even*;

$f(v_{l,j})=f(v_{l-2,j})-q_0-1$,   $j\in[1,n]$, $l\in[1,m+n]$ *l and even*;

$f(u_{l,i})=f(u_{l-2,i})-q_0-1$,   $i\in[1,m]$, $l\in[1,m+n]$ and $l$ *odd*;

$f(v_{l,j})=f(v_{l-2,j})+q_0+1$,   $j\in[1,n]$, $l\in[1,m+n]$ *l and odd.*

Above defined labeling function $f$ give rise edge labels $1,2,\cdots,q_0, q_0+2, q_0+3,\cdots,2q_0+1,2q_0+3,2q_0+4,\cdots,3q_0+2,\cdots,(q_0+1)(m+n)+1,\cdots,q$ to all the copies of $G$. To make $G^\star$ as graceful graph $q_0+1,2(q_0+1),\cdots,(m+n)(q_0+1)$ edge labels to be require.

Now we see that the difference of vertex labels for the central copy $G^{(0)}$ with its other copies $G^{(i)}$, $1\leqslant i\leqslant m+n$ is precisely following sequence.

$$(m+n)(q_0+1)$$
$$q_0+1$$
$$(m+n-1)(q_0+1)$$
$$2(q_0+1)$$
$$\cdots$$
$$(m+n)(q_0+1)/2,\ \text{when } m+n\equiv 0(\text{mod}2)$$
or
$$(m+n-1)(q_0+1)/2,\ \text{when } m+n\equiv 1(\text{mod}2)$$
$$(m+n+1)(q_0+1)/2$$

Using above sequence we can produce required edge labels by joining corresponding vertex of $G^{(0)}$ with its other copies of $G^{(i)}$, $\forall\ i=1,2,\cdots,m+n$.Thus $G^\star$admits a graceful labeling and so it is a graceful graph.

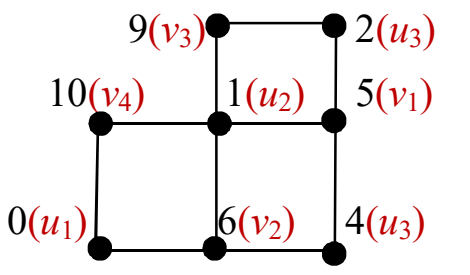

Figure 11 . $G$ and its bipartite graceful labeling

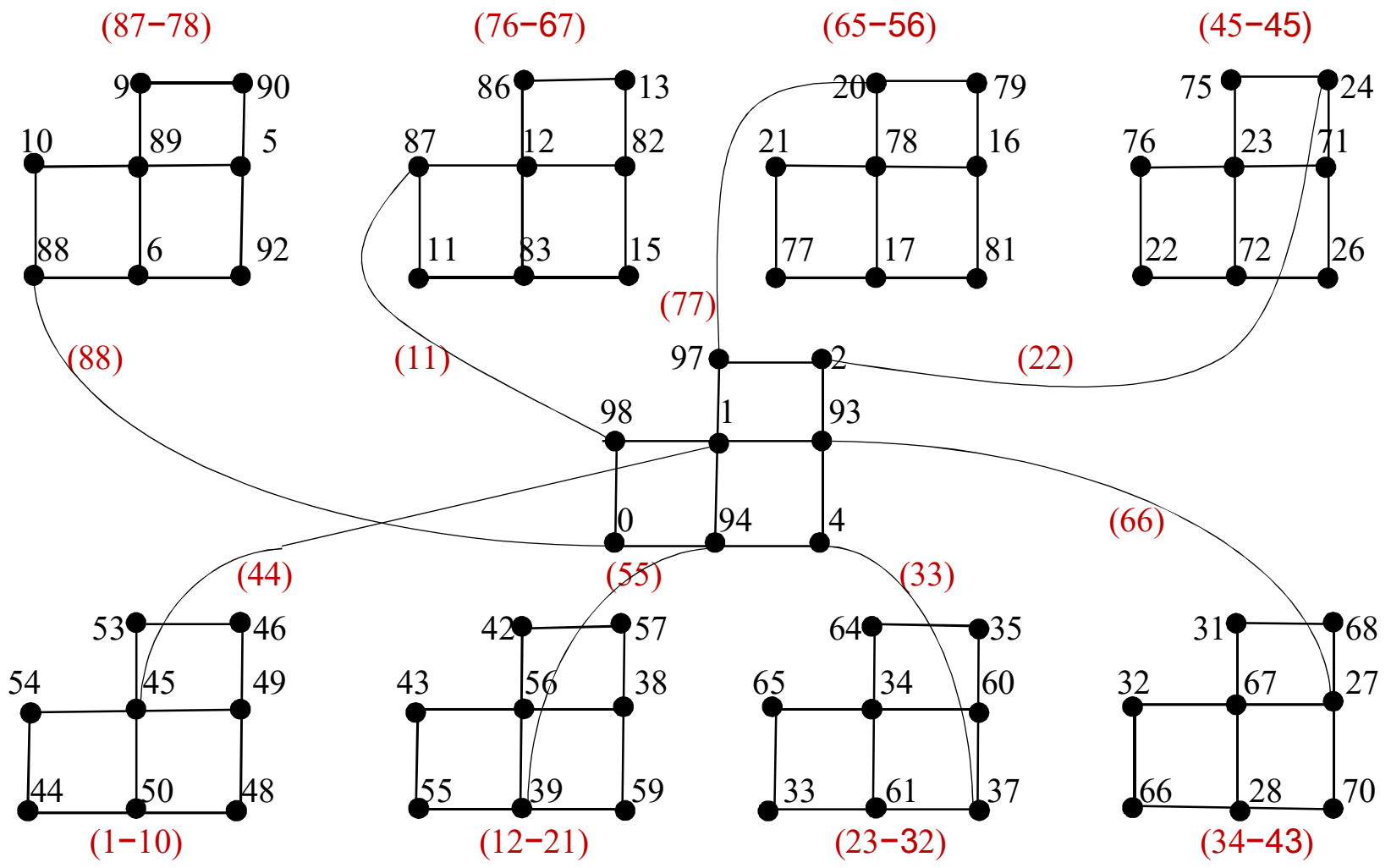

Figure 12 . $G^\star$and its bipartite graceful labeling

By theorem2.5 and lemma 1.1−1.4,respectively,have following consequences.

**Corollary2.5.1** $P_n^\star$of path graph $P_n$ is graceful

**Corollary2.5.2** $C^\star_{4m}$of $C_{4m}$ graphs is graceful graphs.

**Corollary2.5.3**[6] $K^\star_{m,n}$ of complete bipartite graphs is graceful graphs.

**Corollary2.5.4**[12] $(P_r\times P_s)^\star$ of grid graph $P_r\times P_s$) is graceful.